%% file: root.tex
\documentclass[letterpaper, 10 pt, conference]{ieeeconf}  

\usepackage{amsmath, amsfonts, amssymb, arydshln}
\usepackage{algorithm, algorithmic}
\usepackage{latexsym}
\usepackage{dsfont}
\usepackage{floatpag}
\usepackage{float}
\usepackage{tikz}
\usetikzlibrary{shapes,arrows}
\tikzstyle{block} = [draw, rectangle, 
    minimum height=2em, minimum width=4em]
\tikzstyle{sum} = [draw, circle, node distance=1cm]
\tikzstyle{input} = [coordinate]
\tikzstyle{output} = [coordinate]
\tikzstyle{pinstyle} = [pin edge={to-,thin,black}]
\usepackage{verbatim}
\usepackage{cite}
\newtheorem{theorem}{Theorem}
\newtheorem{remark}{Remark}

\usepackage{graphicx,ragged2e}
\DeclareMathOperator*{\argmin}{argmin} 

\IEEEoverridecommandlockouts                              

\title{\LARGE \bf
Closed-Loop Identification of Stabilized Models Using Dual Input-Output Parameterization*}

\author{Ran Chen, Amber Srivastava, Mingzhou Yin, Roy S. Smith 
\thanks{Ran Chen is currently a Master's student at ETH Z\"urich, 8092 Zurich, Switzerland {\tt\small (rachen@student.ethz.ch)}.}
\thanks{Mingzhou Yin, and Roy S. Smith are with the Automatic Control Laboratory (Institut f\"ur Automatik, IfA), ETH Z\"urich, 8092 Zurich, Switzerland {\tt\small(myin/rsmith@control.ee.ethz.ch)}.}%
\thanks{Amber Srivastava is with the Indian Institute of Technology Delhi, 110016 New Delhi, India {(\tt\small asrvstv@mech.iitd.ac.in)}.}
\thanks{*This work was supported by NCCR Automation, funded by the Swiss National Science Foundation (grant number 180545).}
}

\begin{document}

\maketitle
\thispagestyle{empty}
\pagestyle{empty}

\input{sections/00_abs}

\input{sections/01_intro}
\input{sections/02_preli}

\input{sections/03_method}
\input{sections/03.5_review}
\input{sections/04_result}

\input{sections/05_conclu}

\bibliographystyle{IEEEtran}
\bibliography{mybibfile}

\end{document}

%% file: sections/00_abs.tex
\begin{abstract}

This paper introduces a dual input-output parameterization (dual IOP) for the identification of linear time-invariant systems from closed-loop data. It draws inspiration from the recent input-output parameterization developed to synthesize a stabilizing controller. The controller is parameterized in terms of closed-loop transfer functions, from the external disturbances to the input and output of the system, constrained to lie in a given subspace. Analogously, the dual IOP method parameterizes the unknown plant with analogous closed-loop transfer functions, also referred to as dual parameters. In this case, these closed-loop transfer functions are constrained to lie in an affine subspace guaranteeing that the identified plant is \emph{stabilized} by the known controller. Compared with existing closed-loop identification techniques guaranteeing closed-loop stability, such as the dual Youla parameterization, the dual IOP neither requires a doubly-coprime factorization of the controller nor a nominal plant that is stabilized by the controller. The dual IOP does not depend on the order and the state-space realization of the controller either, as in the dual system-level parameterization. Simulation shows that the dual IOP outperforms the existing benchmark methods.
\end{abstract}

%% file: sections/01_intro.tex
\section{Introduction}
System identification, a methodology employed to construct dependable system models from measured data, has widespread applications in the field of engineering \cite{ljung1999system}. This methodology serves as a fundamental step in enabling various aspects of model-based control system design \cite{VANDENHOF1998173}, minimum variance control \cite{GEVERS1986543}, robust control design \cite{zang1995iterative} and more, as outlined in the relevant literature. Open-loop identification methods, also known as direct methods, are often implemented to estimate the system transfer function directly from the input-output data. However, such methods require the plant to be operating and stable in open-loop and may result in an inconsistent estimate if applied to closed-loop data \cite{VANDENHOF1998173, FORSSELL19991215}. Furthermore, if information about the controller is available, it can be leveraged to guarantee the closed-loop stabilizability of the identified plant. This attribute is particularly useful when estimating a plant with a limited stability margin. 

A celebrated indirect method in closed-loop identification is the \textit{Dual Youla Parameterization} (dual YP) \cite{dualyp}. This method requires pre-computing the doubly co-prime factorization of the known controller as well as a nominal plant that is stabilized by the given controller. The dual Youla parameter is estimated as the transfer function from a pair of filtered input and output data. Subsequently, the unknown plant, which is parameterized by the dual Youla parameter, is determined. It is known that if the estimated dual Youla parameter is stable, then the identified plant is stabilized by the known controller. Another indirect method studied recently is the \textit{Dual System-Level Parameterization} (dual SLP) \cite{srivastava2023dual}. In this framework, several closed-loop transfer functions, also known as dual system-level parameters, are defined to parameterize the unknown plant. When constrained to lie in an affine subspace (determined by the state-space realization of the controller), these dual system-level parameters guarantee that the parameterized plant is stabilized by the known controller. As discussed above, the above frameworks have the vital advantage that the identified plant is stabilized by the controller used in the experiment. However, they are dependent on pre-computations (dual YP) and the order of the controller (dual SLS), which may pose additional uncertainties or estimation variance in the identification schemes. 

In fact, the above dual Youla method is motivated by the Youla parameterization (YP) \cite{assadian2022robust} which is extensively applied in robust control, whereas the dual system-level method gains inspiration from the system-level synthesis method (SLS) \cite{ANDERSON2019364, wang2019system} which provides an alternative for parameterized controller design. Recently, the input-output parameterization (IOP) \cite{Furieri_2019} has been developed for controller design, treating the closed-loop transfer functions (from the external disturbances to the input and output of the system) as design variables and exploiting their affine relationships. Specifically, the IOP perspective for controller design is free of any pre-computations such as the doubly-coprime factorization and initial stabilizing controller required in the YP method and the state-space realization of the plant as required in SLS, making it especially effective for large-scale multi-input multi-output (MIMO) systems compared to the others.  Considering the fact that the controller design and the system identification are dual problems, it is worthwhile to investigate the dual application of IOP to system identification.

In this paper, we present the theory of \textit{Dual Input-Output Parameterization} (dual IOP) to obtain a \emph{stabilized} estimate of the plant using closed-loop data. We prove that the set of all plants stabilized by a given controller can be characterized by an affine subspace of four dual input-output parameters. It is shown that identifying these parameters is equivalent to an open-loop identification problem. Similarly, dual IOP does not depend on any pre-computations, hence requiring fewer assumptions than dual YP or dual SLP. In addition, the dual IOP is more computationally efficient than the dual SLP since the former needs fewer optimization variables than the latter. We compare the performance of dual IOP with dual YP and dual SLP via the Bode plots of the estimated plants and the distribution of the estimation errors. Simulation results show that dual IOP exhibits asymptotic convergence in errors and an improvement in estimation, with a decrease in both the median (7\% and 18\% lower than dual SLP and dual YP respectively) and variance of errors. 

It is important to note that there are other well-known closed-loop identification frameworks, including the two-stage method \cite{twostage} and the projection method \cite{forssell2000projection}. We are also aware of several recent developments in closed-loop identification \cite{10.1115/1.4025158, 6981340, AGUERO20111626}. These methods can yield consistent estimates without requiring knowledge of the controller. However, they lack the guarantee that the identified plant is stabilized by the controller. In this work, where the emphasis is on obtaining stabilized models, we exclusively compare our proposed method with relevant literature, namely dual YP and dual SLP.

This paper is structured as follows. Section \ref{sec:preli} provides an overview of the preliminary concepts, including nomenclatures, problem setup, and an introduction to the input-output parameterization. In Section \ref{sec:diop}, we delve into the theory and implementation of dual IOP. This is followed by Section \ref{sec:review}, which offers a comprehensive review of dual SLP and dual YP. Section \ref{sec:result} is dedicated to the comparative analysis of our proposed method against these benchmark closed-loop identification frameworks. Finally, we summarize and draw conclusions in Section \ref{sec:conclu}.

%% file: sections/02_preli.tex
\section{Preliminaries}
\label{sec:preli}
\subsection{Nomenclatures}
We denote vectors and signals by lowercase ($x$) and bold lowercase ($\mathbf{x}$) letters. Uppercase ($A$) and bold uppercase ($\mathbf{A}$) letters represent real matrices and (matrices of) transfer functions, respectively. The operator $z$ denotes a forward shift by one sampling interval as well as the variable in the $z$-domain. We define $\mathcal{RH}_{\infty}$ as the space of all bounded stable and causal transfer functions, and $\mathcal{R}_c$ as the space of all causal transfer functions. Moreover, $\mathcal{RH}_\infty^{p\times m}$ indicates the space of $p\times m$ matrices of the corresponding transfer functions and similarly for $\mathcal{R}_c^{p\times m}$. Finally, $\mathbf{0}$ and $I$ stand for zero matrices and identity matrices with appropriate sizes.

\subsection{Problem Configuration}
\begin{figure}
    \centering
    \begin{tikzpicture}[auto, node distance=2cm,>=latex']
    \node [input, name=input] {};
    \node [sum, right of=input] (sum) {};
    \node [block, right of=sum,
            node distance=2.5cm] (G) {$\mathbf{G}(z)$};
    \draw (sum.west) -- (sum.east)
	(sum.north) -- (sum.south);
    \draw [->] (sum) -- node[name=u] {$u(k)$} (G);
    \node [sum, right of=G, node distance=1.5cm] (disturbance) {};
    \draw (disturbance.west) -- (disturbance.east)
	(disturbance.north) -- (disturbance.south);
    \node [block, above of=disturbance, node distance=1.5cm] (H) {$\mathbf{H}(z)$};
    \node [input, left of=H] (dist_input){$e(k)$};
    \node [output, right of=disturbance] (output) {};
    \node [block, below of=G, node distance=1cm] (measurements) {$\mathbf{K}(z)$};

    \draw [draw,->] (input) -- node {$r(k)$} (sum);
    \draw [->] (G) -- (disturbance);
    \draw [->] (disturbance) -- node [name=y] {$y(k)$}(output);
    \draw [->] (H) -- node {$v(k)$} (disturbance);
    \draw [draw, ->] (dist_input) -- node {$e(k)$} (H);
    \draw [->] (y) |- (measurements);
    \draw [->] (measurements) -| (sum);
\end{tikzpicture}
    \caption{Block diagram of an LTI feedback system.}
    \label{fig:block}
\end{figure}
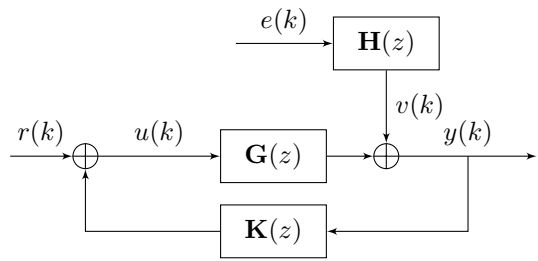

We consider a closed-loop linear time-invariant (LTI) system in discrete time given by
\begin{equation}
    \label{eq:iodynamics} 
        \begin{aligned}
        y(k) & = \mathbf{G}(z)u(k) + v(k), \\
        u(k)  & = \mathbf{K}(z)y(k) + r(k),
    \end{aligned}    
\end{equation}
which is incorporated with a positive linear feedback controller $\mathbf{K}$. The controller $\mathbf{K}(z)\in\mathcal{R}_c$ is exactly known, whereas the input-output (IO) transfer function $\mathbf{G}(z)\in \mathcal{R}_c$ is unknown and to be identified. A block diagram of this system is given in Figure \ref{fig:block}. Note $\mathbf{G}$ and $\mathbf{K}$\footnote{For conciseness, we omit the shift index $z$ for all transfer functions in the following.} are not necessarily internally stable, but $\mathbf{K}$ internally stabilizes $\mathbf{G}$. The output $y(k)\in\mathbb{R}^p$, control input $u(k)\in\mathbb{R}^m$, and reference $r(k)\in\mathbb{R}^m$ are accessible at each time $k=0,\cdots, N-1$. The output is further corrupted by the unmeasurable signal $v(k)=\mathbf{H}e(k)$ with independent and identically distributed noise $e(k)$.

\subsection{Input-Output Parameterization}
This subsection reviews the principles of the IOP method used for output feedback controller design \cite{Furieri_2019}, as it is fundamental to our proposed dual IOP framework. In the IOP framework, the controller $\mathbf{K}$ is to be designed whereas the plant $\mathbf{G}$ is known. Considering the external signal pair $(\mathbf{v}, \mathbf{r}$) and the dependent signal pair $(\mathbf{y}, \mathbf{u}$), (\ref{eq:iodynamics}) leads to the following closed-loop equations,
\begin{equation}
    \begin{aligned}
    \begin{bmatrix}
        \mathbf{y} \\
        \mathbf{u}  
    \end{bmatrix}
    & =
    \begin{bmatrix}
         (I-\mathbf{GK})^{-1}  & (I-\mathbf{GK})^{-1}\mathbf{G} \\
         (I-\mathbf{KG})^{-1}\mathbf{K} & (I-\mathbf{KG})^{-1} 
    \end{bmatrix}
    \begin{bmatrix}
        \mathbf{v} \\
        \mathbf{r}   
    \end{bmatrix} \\ 
    &=
    \begin{bmatrix}
        \mathbf{W} & \mathbf{X} \\
        \mathbf{Y} & \mathbf{Z}
    \end{bmatrix}
    \begin{bmatrix}
        \mathbf{v} \\
        \mathbf{r}  
    \end{bmatrix},
    \end{aligned}
\label{eq:XWYZ}
\end{equation}
where $\mathbf{W}, \mathbf{X}, \mathbf{Y}, \mathbf{Z}$ define the parameterization of the four closed-loop transfer functions in terms of the system $(\mathbf{G}, \mathbf{K})$. Under the assumption that $\mathbf{K}$ internally stabilizes $\mathbf{G}$, these closed-loop transfer functions are stable and causal \cite{francis1987course}, i.e. $\mathbf{W}, \mathbf{X}, \mathbf{Y}, \mathbf{Z} \in \mathcal{RH}_\infty$. In the IOP framework,  these transfer functions are treated as optimization variables, before using their relationship to indirectly determine $\mathbf{K}=\mathbf{YW}^{-1}$. Let the set of internally stabilizing controllers be given by,
\begin{equation}
    \mathcal{C}_{\mathrm{stab}}(\mathbf{G}) = \left\{ \mathbf{K}\in \mathcal{R}_c \;|\; \mathbf{K}\;\mathrm{stablizes\;\mathbf{G}}  \right\},
\end{equation}
then following theorem summarizes the properties of IOP.

\begin{theorem}
    Consider the LTI system in (\ref{eq:iodynamics}), the following statements are true \cite{Furieri_2019}.

    \begin{enumerate}
        \item For any arbitrary controller $\mathbf{K}\in\mathcal{C}_{\mathrm{stab}}(\mathbf{G}) $, there always exist $(\mathbf{W}, \mathbf{X}, \mathbf{Y}, \mathbf{Z})$ that lie in the following affine subspace,
        \begin{subequations}
        \label{eq:affine}
            \begin{align}
                \begin{bmatrix}
                I & -\mathbf{G}
                \end{bmatrix} 
                \begin{bmatrix}
                    \mathbf{W} & \mathbf{X} \\
                    \mathbf{Y} & \mathbf{Z}
                \end{bmatrix}  & = 
                \begin{bmatrix}
                    I & \mathbf{0}
                \end{bmatrix},\\ 
                \begin{bmatrix}
                    \mathbf{W} & \mathbf{X} \\
                    \mathbf{Y} & \mathbf{Z}
                \end{bmatrix} 
                \begin{bmatrix}
                -\mathbf{G} \\ I
                \end{bmatrix}  & = 
                \begin{bmatrix}
                     \mathbf{0} \\ I
                \end{bmatrix}, \\
                 \mathbf{W}, \mathbf{X}, \mathbf{Y}, \mathbf{Z}& \in\mathcal{RH}_\infty.
            \end{align}
        \end{subequations}

        \item For any transfer functions $(\mathbf{W}, \mathbf{X}, \mathbf{Y}, \mathbf{Z})$ that lie in the
        affine subspace (\ref{eq:affine}a)-(\ref{eq:affine}c), the parameterized controller $\mathbf{K}=\mathbf{YW}^{-1}=\mathbf{Z}^{-1}\mathbf{Y}\in\mathcal{C}_{\mathrm{stab}}(\mathbf{G}) $, i.e. the closed-loop system (\ref{eq:iodynamics}) is stable.
    \end{enumerate}
\end{theorem}

Therefore, the controller $\mathbf{K}$ can be determined by finding the closed-loop transfer functions that minimize a pre-defined cost function (see \cite{Furieri_2019} and \cite{zheng2022system}) subject to subspace constraints (\ref{eq:affine}a)-(\ref{eq:affine}c). If the optimal solution $(\hat{\mathbf{W}}, \hat{\mathbf{X}}, \hat{\mathbf{Y}}, \hat{\mathbf{Z}})$ exists, the resulting $\hat{\mathbf{K}}$ is guaranteed to stabilize $\mathbf{G}$. Hence, the controlled system has closed-loop stability. In the next section, we exploit the duality of controller synthesis and system identification, and propose the novel dual input-output parameterization (dual IOP) which chooses the optimal estimate $\hat{\mathbf{G}}$ from the set of all plants that are stabilized by a given output feedback controller $\mathbf{K}$.

%% file: sections/03_method.tex
\section{The Dual Input-Output Parameterization}
\label{sec:diop}
\subsection{Theory of Dual IOP}
In the dual IOP framework, we reverse the role of the plant and the controller, i.e., we look for a plant model that \emph{is stabilized by} the given controller. To see this, we consider the following dual system labeled by a tilde: $ (\tilde{\mathbf{G}}, \tilde{\mathbf{K}}) = (\mathbf{K},\mathbf{G})$, $(\tilde{\mathbf{v}}, \tilde{\mathbf{r}}) = (\mathbf{r},\mathbf{v})$ and $(\tilde{\mathbf{y}}, \tilde{\mathbf{u}}) = (\mathbf{u},\mathbf{y})$, thus the objective is to {\em synthesize} $\tilde{\mathbf{K}}$ (i.e., estimate $\mathbf{G}$) such that the closed-loop system $(\tilde{\mathbf{G}}, \tilde{\mathbf{K}})$ is stable. Substituting the reversed system into (\ref{eq:XWYZ}) results in
\begin{equation}
\label{eq:WXYZ2}
\begin{aligned}
    \begin{aligned}
    \begin{bmatrix}
        \tilde{\mathbf{y}} \\
        \tilde{\mathbf{u}}
    \end{bmatrix}
    & = 
    \begin{bmatrix}
         (I-\tilde{\mathbf{G}}\tilde{\mathbf{K}})^{-1}  & (I-\tilde{\mathbf{G}}\tilde{\mathbf{K}})^{-1}\tilde{\mathbf{G}} \\
         (I-\tilde{\mathbf{K}}\tilde{\mathbf{G}})^{-1}\tilde{\mathbf{K}} & (I-\tilde{\mathbf{K}}\tilde{\mathbf{G}})^{-1}
    \end{bmatrix} 
    \begin{bmatrix}
        \tilde{\mathbf{v}} \\
        \tilde{\mathbf{r}}   
    \end{bmatrix} \\
    \Rightarrow \begin{bmatrix}
    \mathbf{u}\\
    \mathbf{y}
        \end{bmatrix}
    &=
    \begin{bmatrix}
          (I-\mathbf{KG})^{-1} & (I-\mathbf{KG})^{-1}\mathbf{K}  \\
          (I-\mathbf{GK})^{-1}\mathbf{G} & (I-\mathbf{GK})^{-1} 
    \end{bmatrix}
    \begin{bmatrix}
        \mathbf{r} \\
        \mathbf{v}   
    \end{bmatrix} \\ 
    &=
    \begin{bmatrix}
         \mathbf{Z} & \mathbf{Y} \\
            \mathbf{X} & \mathbf{W}
    \end{bmatrix}
    \begin{bmatrix}
        \mathbf{r} \\
        \mathbf{v} 
    \end{bmatrix}.
    \end{aligned} \\
\end{aligned}
\end{equation}

Considering $\tilde{\mathbf{G}}=\mathbf{K}$, the constraints in (\ref{eq:affine}) become
\begin{subequations}
\label{eq:affine3}
    \begin{align}
        \begin{bmatrix}
        I & -\mathbf{K}
        \end{bmatrix} 
        \begin{bmatrix}
            \mathbf{Z} & \mathbf{Y} \\
            \mathbf{X} & \mathbf{W}
        \end{bmatrix}  & = 
        \begin{bmatrix}
            I & 0
        \end{bmatrix}\\ 
        \begin{bmatrix}
            \mathbf{Z} & \mathbf{Y} \\
            \mathbf{X} & \mathbf{W}
        \end{bmatrix} 
        \begin{bmatrix}
        -\mathbf{K} \\ I
        \end{bmatrix}  & = 
        \begin{bmatrix}
             0 \\ I
        \end{bmatrix} \\
        \mathbf{Z}, \mathbf{Y}, \mathbf{X}, \mathbf{W}& \in\mathcal{RH}_\infty.
    \end{align}
\end{subequations}

Thus, given $\tilde{\mathbf{K}}=\mathbf{G}=\mathbf{W}^{-1}\mathbf{X}$ with $\mathbf{W}, \mathbf{X}$ satisfying (\ref{eq:affine3}), the closed-loop system $(\tilde{\mathbf{G}},\tilde{\mathbf{K}})$ (equivalently, $(\mathbf{K}, \mathbf{G})$) is stable. Thus, similar to Theorem 1, the idea behind dual IOP is summarized in the following theorem.
\begin{theorem}
\label{thm:thm2}
    Consider the LTI system in (\ref{eq:iodynamics}), the following statements are true.
    \begin{enumerate}
        \item Suppose a controller $\mathbf{K}$ is given and stabilizes the unknown plant $\mathbf{G}$. Then there always exist $(\mathbf{W}, \mathbf{X}, \mathbf{Y}, \mathbf{Z})$ that lie in the following affine subspace,
        \begin{subequations}
        \label{eq:affine_diop}
            \begin{align}
                \begin{bmatrix}
                -\mathbf{K} & I
                \end{bmatrix} 
                \begin{bmatrix}
                    \mathbf{W} & \mathbf{X} \\
                    \mathbf{Y} & \mathbf{Z}
                \end{bmatrix}  & = 
                \begin{bmatrix}
                    0 & I
                \end{bmatrix},\\ 
                \begin{bmatrix}
                    \mathbf{W} & \mathbf{X} \\
                    \mathbf{Y} & \mathbf{Z}
                \end{bmatrix} 
                \begin{bmatrix}
                I \\ -\mathbf{K}
                \end{bmatrix}  & = 
                \begin{bmatrix}
                     I \\ 0
                \end{bmatrix}, \\
                 \mathbf{W}, \mathbf{X}, \mathbf{Y}, \mathbf{Z}& \in\mathcal{RH}_\infty.
            \end{align}
        \end{subequations}

        \item For any transfer functions $(\mathbf{W}, \mathbf{X}, \mathbf{Y}, \mathbf{Z})$ that lie in the
        affine subspace (\ref{eq:affine_diop}a)-(\ref{eq:affine_diop}c), the parameterized plant $\mathbf{G}=\mathbf{W}^{-1}\mathbf{X} = \mathbf{XZ}^{-1}\in\mathcal{G}_{\mathrm{stab}}(\mathbf{K}) $, where $\mathcal{G}_{\mathrm{stab}}(\mathbf{K})$ is the set of all plants stabilized by $\mathbf{K}$.
    \end{enumerate}
\end{theorem}

The proof of Theorem 2 is analogous to that of Theorem 1 in \cite{Furieri_2019}, thus omitted here. Subsequently, we formulate the dual IOP closed-loop identification framework by posing the following optimization problem,
\begin{equation}
\label{prob:d-iop}
    \begin{aligned}
        \min_{\mathbf{W}, \mathbf{X}, \mathbf{Y}, \mathbf{Z}}\; & f(\mathbf{r}, \mathbf{y}, \mathbf{X})  \\
        \mathrm{s.t.}\; & (\ref{eq:affine_diop}a)-(\ref{eq:affine_diop}c),
    \end{aligned}
\end{equation}
where $f(\mathbf{r}, \mathbf{y}, \mathbf{X})$ is a function representing the model fitting error between $\mathbf{r}$ and $\mathbf{y}$. A typical choice of $f$ (in the absence of the knowledge of the noise filter $\mathbf{H}$) is the two-norm residuals, i.e., $f=\left\|\mathbf{y} - \mathbf{X}\mathbf{r} \right\|_2^2=\left\|\boldsymbol{\epsilon}\right\|_2^2$, where $\mathbf{\epsilon}$ is the prediction error, as identical to the convex cost function defined in \cite{srivastava2023dual}. Clearly, the optimization problem (\ref{prob:d-iop}) translates the identification of the plant $\mathbf{G}$ into an equivalent open-loop identification problem of $\mathbf{W}, \mathbf{X}, \mathbf{Y}, \mathbf{Z}$ with affine constraints, where the external reference $\mathbf{r}$ and the measured output $\mathbf{y}$ are statistically uncorrelated. Hence, the dual input-output parameters are estimated via the convex program in (\ref{prob:d-iop}).
\begin{remark}
    It is assumed that $\mathbf{r}$ and $\mathbf{e}$ (thus $\mathbf{v}=\mathbf{He}$) are uncorrelated in (\ref{eq:iodynamics}). Suppose $\mathcal{W}, \mathcal{X}, \mathcal{Y}, \mathcal{Z}$ represent the model classes that contain all admissible solutions of (\ref{prob:d-iop}), and suppose the true dual input-output parameters belong to these classes. Then, $\mathbf{W}, \mathbf{X}, \mathbf{Y}, \mathbf{Z}$ can be consistently and unbiasedly \cite{FORSSELL19991215} estimated by (\ref{prob:d-iop}), provided $\mathbf{r}$ is persistently exciting.
\end{remark}
\begin{remark}
     An implicit advantage of the dual IOP framework with respect to its competitors is that the dual IOP directly parameterizes the closed-loop transfer functions from noises and references to inputs and outputs. Any knowledge about the closed-loop behaviors of the original system can be converted into equivalent constraints and imposed on the dual input-output parameters. One of the examples is the identification of systems that are positively stabilized by the given controller \cite{rantzer2015scalable}, which means the resulting closed-loop transfer functions (e.g. $\mathbf{W,X,Y,Z}$) are positive. This prior information can be easily incorporated as constraints on the dual parameters in the optimization problem (\ref{prob:d-iop}). Nevertheless, the investigation of identifying systems with specific closed-loop properties using the dual IOP is out of the scope of this paper. It remains a potential direction for future work.
\end{remark}

\subsection{Implementation of Dual IOP}
A straightforward implementation of the dual IOP identification framework is introduced as follows. It is worth noting that the optimization problem (\ref{prob:d-iop}) is generally infinite-dimensional. To avoid solving such an infinite-dimensional problem, we consider a length $\tau $ finite impulse response (FIR) parameterization to approximate the closed-loop transfer functions $\mathbf{W}, \mathbf{X}, \mathbf{Y}, \mathbf{Z}$, namely, 
\begin{equation}
\label{eq:fir}
    \begin{aligned}
        \mathbf{W}=\sum_{k=0}^{\tau -1}W[k]z^{-k}, & \;\;\mathbf{X}=\sum_{k=0}^{\tau -1}X[k]z^{-k}, \\
        \mathbf{Y}=\sum_{k=0}^{\tau -1}Y[k]z^{-k}, & \;\;\mathbf{Z}=\sum_{k=0}^{\tau -1}Z[k]^{-k}, 
    \end{aligned}
\end{equation}
where $W[k]$, $X[k]$, $Y[k]$, $Z[k]$ are real matrices serving as decision variables which represent the $k$-th element of the corresponding transfer function $\forall\;k=0,\cdots,\tau-1$. 

In the discrete-time case, it has been proved that for $\tau \rightarrow \infty$, the FIR parameterization in (\ref{eq:fir}) spans the entire $\mathcal{R}_c$ \cite{pohl2009advanced}, and the same parameterization as above is exploited in \cite{Furieri_2019} as well for the case of controller design. The approximation order $\tau$ is then a hyperparameter that is to be chosen appropriately: an excessively large $\tau$ leads to huge computational complexity and undesired overfitting of noise, whereas $\tau$ being too short makes the truncation error large, thus deteriorating the quality of the estimate. 

Subsequently, the stability and causality constraint (\ref{eq:affine_diop}c) is inherently satisfied since a transfer function in FIR parameterization has poles only at the origin. The affine subspace constraints (\ref{eq:affine_diop}a)-(\ref{eq:affine_diop}b) can be implemented as a group of linear equality constraints. An example for $\mathbf{Z}-\mathbf{KX}=I $ is given as 
\begin{equation}
\label{eq:Z-KX=I}
    \begin{bmatrix}
        Z[0] \\
        Z[1] \\
        \vdots \\
        Z[\tau - 1]\\
        0\\
        \vdots\\
        0
    \end{bmatrix} - 
    \mathrm{Co}(\mathbf{K})\cdot 
    \begin{bmatrix}
        X[0] \\
        X[1] \\
        \vdots \\
        X[\tau - 1]
    \end{bmatrix}=
    \begin{bmatrix}
        I \\
        \mathbf{0}
    \end{bmatrix},
\end{equation}
where $\mathrm{Co}(\mathbf{K})\in \mathbb{R}^{m\cdot(\nu + \tau -1)\times p\cdot \tau }$ represents the convolution matrix of $\mathbf{K}$ and is constructed by,
\begin{equation}
    \mathrm{Co}(\mathbf{K}) = \begin{bmatrix}
        K[0] & 0  & \cdots  & 0\\
        K[1] & K[0]  & \cdots & 0  \\
          \vdots     &    \vdots  &    \ddots & \vdots \\
        K[\nu -1] & K[\nu -2] & \cdots & K[0] \\
        0 & K[\nu -1] & \cdots & K[1] \\
        \vdots & \vdots & \ddots & \vdots \\
        0 & 0 & \cdots & K[\nu -1 ] 
    \end{bmatrix}.
\end{equation}

The equation above is valid when $\mathbf{K}$ is a deadbeat controller, e.g. $\mathbf{K}(z) = \sum_{k=0}^{\nu -1}K[k]z^{-k}$ for some finite order $\nu \leq \tau$. In case $\mathbf{K}$ has a general fractional structure, i.e. $\mathbf{K} = \mathbf{D}_K^{-1}\mathbf{N}_K $, we have $\mathbf{Z}-\mathbf{KX} \Rightarrow \mathbf{D}_K\mathbf{Z}-\mathbf{N}_K\mathbf{X}=\mathbf{D}_K $, which also results in a group of linear equality constraints. The cost function $f=\left\|\mathbf{y} - \mathbf{X}\mathbf{r} \right\|_2^2$ is constructed as,
\begin{equation}
\label{eq:cost}
    f=\left\| \begin{bmatrix}
        y(0)^\top \\
        y(1)^\top  \\
        \vdots \\
        y(N-1)^\top 
    \end{bmatrix}  - \mathrm{toep}(\mathbf{r}^\top)\cdot
    \begin{bmatrix}
        X[0]^\top \\
        X[1]^\top \\
        \vdots \\
        X[\tau -1]^\top \\
    \end{bmatrix}\right\|_2^2,
\end{equation}
where $\mathrm{toep}(\mathbf{r^\top})\in\mathbb{R}^{N\times {m\cdot \tau}}$ represents the (block) Toeplitz matrix of the transpose of $\mathbf{r}$ with the first column block being $[r(0),r(1),\cdots,r(N-1)]^\top$ and the first row block being $[r(0)^\top,\mathbf{0}]$.
The implementation of the dual IOP framework is summarized by Algorithm 1.
\begin{algorithm}
    \caption{}
 \begin{algorithmic}[1]
 \renewcommand{\algorithmicrequire}{\textbf{Input:}}
 \renewcommand{\algorithmicensure}{\textbf{Output:}}
 \REQUIRE $\mathbf{y}$, $\mathbf{r}$, $\mathbf{K}$, $\tau$
 \ENSURE  $\hat{\mathbf{G}}$
  \STATE Determine input (output) number $m$ ($p$) from $\mathbf{K}$
  \STATE Define matrix variables $W[k]\in\mathbb{R}^{p\times p}$, $X[k]\in\mathbb{R}^{p\times m}$, $Y[k]\in\mathbb{R}^{m\times p}$, $Z[k]\in\mathbb{R}^{m\times m}$ $\forall k=0,\cdots,\tau -1$
 \STATE Construct the cost function $f$ via (\ref{eq:cost})
 \STATE Construct the affine constraints in (\ref{eq:affine_diop}a), (\ref{eq:affine_diop}b) via (\ref{eq:Z-KX=I})
 \STATE Solve the quadratic program with the cost and constraints defined in Steps 3 and 4
 and obtain the estimated FIR elements $\hat{W}[k],\hat{X}[k], \hat{Y}[k], \hat{Z}[k]$
 \STATE Construct $\hat{\mathbf{W}}=\sum_{k=0}^{\tau -1}\hat{W}[k]z^{-k}$ and similar for $\hat{\mathbf{X}}$
 \STATE Compute $\hat{\mathbf{G}} = \hat{\mathbf{W}}^{-1}\hat{\mathbf{X}}$ and return $\hat{\mathbf{G}}$
 \end{algorithmic} 
\end{algorithm}

%% file: sections/03.5_review.tex
\section{Review of the Benchmark Methods}
\label{sec:review}
This section provides a brief overview of two existing closed-loop identification methods which also guarantee the identified plant is stabilized by the known controller, namely the dual YP and dual SLP methods. 

\subsection{The Dual YP Identification}
Given an initial guess of the plant $\mathbf{G}_X$ stabilized by $\mathbf{K}$, suppose $\mathbf{K}$ and $\mathbf{G}_X$ admit coprime factorizations $\mathbf{K}=\mathbf{D}_K^{-1}\mathbf{N}_K$ and $\mathbf{G}_X=\mathbf{D}_X^{-1}\mathbf{N}_X$, respectively \cite{copfac}. Subsequently, we construct the virtual input $\boldsymbol{\alpha}$ and output $\boldsymbol{\beta}$ by,
\begin{equation}
    \begin{aligned}
        \boldsymbol{\beta} & = \mathbf{D}_X \mathbf{y} - \mathbf{N}_X \mathbf{u},\\
        \boldsymbol{\alpha} & = \mathbf{D}_K \mathbf{u} - \mathbf{N}_K \mathbf{y}= \mathbf{D}_K\mathbf{r}.
    \end{aligned}
\end{equation}

Since $\boldsymbol{\alpha}$ only depends on the reference $\mathbf{r}$, it is statistically uncorrelated with the noise $\mathbf{v}$. Hence, the transfer function $\mathbf{Q}$ from $\boldsymbol{\alpha}$ to $\boldsymbol{\beta}$, also known as the dual Youla parameter, can be consistently identified by the following open-loop equation,
\begin{equation}
\label{eq:d-yp}
    \boldsymbol{\beta} = \mathbf{Q}\boldsymbol{\alpha} + \mathbf{Se},
\end{equation}
where $\mathbf{S}$ represents the noise dynamics of this equivalent open-loop system. The dual Youla parameter $\mathbf{Q}$ can be identified by least squares via,
\begin{equation}
    \hat{\mathbf{Q}} = \argmin_{\mathbf{Q}\in\mathcal{RH}_\infty}\;\|\boldsymbol{\beta}-\mathbf{Q}\boldsymbol{\alpha}\|_2,
\end{equation}
and the estimated plant $\hat{\mathbf{G}}$ is given by,
\begin{equation}
    \hat{\mathbf{G}} = \left( \mathbf{D}_X + \hat{\mathbf{Q}}\mathbf{N}_K \right)^{-1}\left( \mathbf{N}_X + \hat{\mathbf{Q}}\mathbf{D}_K \right),
\end{equation}
which is stabilized by $\mathbf{K}$ if and only if $\hat{\mathbf{Q}}$ is stable.

\subsection{The Dual SLP Identification}
\label{sec:dslp}
The dual SLP framework starts by realizing a state-space representation $(A_K, B_{K}, C_K, D_{K})$ of $\mathbf{K}$ so that the dynamics of the controller can be described by,
\begin{equation}
\begin{aligned}
        \xi[k+1] & = A_K \xi(k) + B_{K} y(k) + B_{Kv}v(k), \\
        u(k) &= C_K\xi(k) + D_{K}y(k) + D_{Kv}v(k),
\end{aligned}
\end{equation}
where $\xi(k)$ is the internal state and $B_{Kv}$,$D_{Kv}$ represent the noise dynamics exerting on $\xi(k)$ and the control input $u(k)$, respectively. Consequently, the dual SLP identification framework solves the following optimization problem,
\begin{subequations}
\label{eq:dslp}
\begin{align}
    & \min_{\mathbf{R}_K,\mathbf{N}_K,\mathbf{M}_K,\mathbf{L}_K} \;\; f(\mathbf{r},\mathbf{y},\mathbf{L}_K) \\
    \mathrm{s.t.}& \; 
    \begin{bmatrix}
        zI-A_K  & -B_K
    \end{bmatrix}
    \begin{bmatrix}
        \mathbf{R}_K & \mathbf{N}_K \\
        \mathbf{M}_K & \mathbf{L}_K
    \end{bmatrix} = 
    \begin{bmatrix}
        I & 0
    \end{bmatrix},\\
    & \;
    \begin{bmatrix}
        \mathbf{R}_K & \mathbf{N}_K \\
        \mathbf{M}_K & \mathbf{L}_K
    \end{bmatrix} \begin{bmatrix}
        zI-A_K  \\ -C_K
    \end{bmatrix} = 
    \begin{bmatrix}
        I \\ 0
    \end{bmatrix},\\
    & \mathbf{R}_K,\mathbf{N}_K,\mathbf{M}_K\in\frac{1}{z}\mathcal{RH}_\infty,\; \mathbf{L}_K\in\mathcal{RH}_\infty.
    \end{align}
\end{subequations}

Specifically, the transfer functions $\mathbf{R}_K\in\frac{1}{z}\mathcal{RH}_\infty^{n\times n}$, $\mathbf{N}_K\in\frac{1}{z}\mathcal{RH}_\infty^{n\times m}$, $\mathbf{M}_K\in\frac{1}{z}\mathcal{RH}_\infty^{p\times n}$, $\mathbf{L}_K\in\mathcal{RH}_\infty^{p\times m}$ are called the closed-loop response functions. Moreover, $\mathbf{L}_K$ is proved to be the transfer function from $\mathbf{r}$ to $\mathbf{y}$, i.e.,
\begin{equation}
    \mathbf{y} = (I-\mathbf{GK})^{-1}\mathbf{G}\mathbf{r} = \mathbf{L}_K\mathbf{r}.
\end{equation}

Therefore, the dual SLP parameter $\mathbf{L}_K$ is identical to the dual IOP parameter $\mathbf{X}$, indicating the relation between these two parameterizations. Supposing that the optimal solution of Problem \ref{eq:dslp} is $(\hat{\mathbf{R}}_K, \hat{\mathbf{N}}_K, \hat{\mathbf{M}}_K, \hat{\mathbf{L}}_K)$, the estimated plant can be constructed as,
\begin{equation}
\label{eq:Ghatslp}
    \hat{\mathbf{G}} = \hat{\mathbf{L}}_K - \hat{\mathbf{M}}_K \hat{\mathbf{R}}_K^{-1} \hat{\mathbf{N}}_K.
\end{equation}

\begin{remark}
     As proved in\cite{srivastava2023dual}, $\hat{\mathbf{L}}_K$ and the optimization problem (\ref{eq:dslp}) remains invariant to the state-space realization of $\mathbf{K}$. However, if the controller has a large order, i.e., the dimensions of the transition matrices $A_K, B_K,C_K$ are large, then the dimensions of the dual parameters will be large too. Thus, the size of the decision variable space (or, the computational complexity) in dual SLS depends on the order of the given controller $\mathbf{K}$. More details are elaborated in Section \ref{sec:simu}.

\end{remark}

%% file: sections/04_result.tex
\section{Numerical Results}
\label{sec:result}
In this section, the performance of the proposed dual IOP identification framework is illustrated using simulation results and compared with dual YP and dual SLP.

\subsection{Simulation Configuration}
\label{sec:simu}
We consider the following single-input single-output (SISO) system \cite{srivastava2023dual,twostage} given by
\begin{equation}
\label{eq:emp_sys}
\begin{aligned}
    \mathbf{G}_0 & = \frac{1}{1-1.6z^{-1}+0.89z^{-2}}, \\
    \mathbf{K}_0 & = z^{-1}-0.8z^{-2}, \\
    \mathbf{H}_0 & = \frac{1-1.56z^{-1} + 1.045z^{-2}-0.3338}{1- 2.35z^{-1} + 2.09z^{-2} - 0.6675}. \\
\end{aligned}
\end{equation}
The estimation error to quantify the quality of the estimate $\hat{\mathbf{G}}$ of the plant $\mathbf{G}_0$ is defined as
\begin{equation}
\label{eq:err}
\begin{aligned}
    \mathrm{err}(\hat{\mathbf{G}}) = \sum_{i=1}^{L}\frac{\|\mathbf{G}_0(j\omega_i) - \hat{\mathbf{G}}(j\omega_i)\|_2}{\|\mathbf{G}_0(j\omega_i)\|_2} \times 100\%,
\end{aligned}
\end{equation}
where $N$ is the data length and $\omega_i$'s are $L$ equally spaced positive frequencies within $(0,\pi)$ with $L=(N+1)/2$. We further choose the pseudorandom binary sequence (PRBS) \cite{naszodi1979digital} with a magnitude of 1 as the reference excitation $r(k)$. For PRBS signals, $N$ is restricted to $2^d-1$, $\forall d\in \mathbb{N}$, and $d$ is chosen from 8 to 14. The noise sequence $e(k)$ is sampled from the normal distribution $\mathcal{N}(0,1)$, $\forall k=0,\cdots, N-1$. We perform the experiments 100 times (i.e., with 100 different realizations of the noise) to illustrate the robustness of the dual IOP methods as will illustrated in Figure \ref{fig:errcomp}. 

The model structures of the dual IOP parameters are chosen as $\tau$-ordered FIR as mentioned in (\ref{eq:fir}) with $\tau= 14$. Similarly, the dual SLP parameters $\mathbf{R}_K, \mathbf{N}_K, \mathbf{M}_K, \mathbf{L}_K$ are also modeled as FIR functions,
\begin{equation}
    \begin{aligned}
        \mathbf{R}_K=\sum_{k=1}^{\tau -1}R_K[k]z^{-k},\;\; & \mathbf{N}_K=\sum_{k=1}^{\tau -1}N_K[k]z^{-k}, \\
        \mathbf{M}_K=\sum_{k=1}^{\tau -1}M_K[k]z^{-k},\; \;& \mathbf{L}_K=\sum_{k=0}^{\tau -1}L_K[k]z^{-k},
    \end{aligned}
\end{equation}
using the same $\tau$. Note that, $\mathbf{R}_K$, $\mathbf{N}_K$, $\mathbf{M}_K$ do not have the zeroth element since they belong to $\frac{1}{z}\mathcal{RH}_\infty$ as shown in (\ref{eq:dslp}d). Similarly, the dual Youla parameter of dual YP is parameterized as
\begin{equation}
    \mathbf{Q}=\sum_{k=0}^{\tau -1}Q[k]z^{-k},
\end{equation}
with $\tau=14$. The following transfer functions are chosen as the initial nominal plant $\mathbf{G}_X$ in dual YP to test its performance under different scenarios:
\begin{itemize}
    \item $\mathbf{G}_a$ = 0, a zero gain.
    \item $\mathbf{G}_b=\hat{\mathbf{G}}_{\mathrm{DY}}$, the estimated plant by dual YP using additional data and $\mathbf{G}_a$ above as the initial plant. This choice is equivalent to a scenario where certain prior knowledge of the plant is provided.
    \item $\mathbf{G}_c = -\frac{1}{z+0.5}$, an arbitrarily selected function that is stabilized by $\mathbf{K}$ as used in \cite{srivastava2023dual}.
\end{itemize}

\subsection{Results and Discussions}
\label{sec:discussion}
\begin{figure}
    \centering
    \includegraphics[width = 0.49\textwidth]{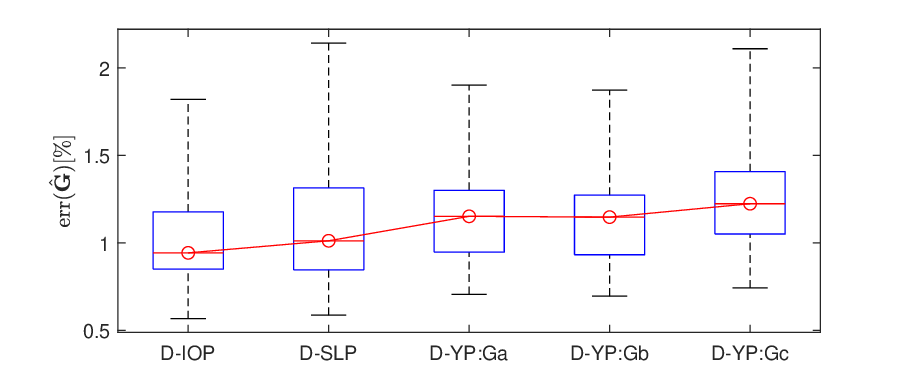}
    \caption{Box plots of $\mathrm{err}(\hat{\mathbf{G}})$ distributions in [\%] using dual IOP, dual SLP, dual YP with $\mathbf{Gx}=\mathbf{G}_a$, $\mathbf{G}_b$ and $\mathbf{G}_c$ separately. }
    \label{fig:errcomp}
\end{figure}

We collect the estimated plants $\hat{\mathbf{G}}$ from all three methods using independently generated data, calculate the errors as defined in (\ref{eq:err}), and visualize the results in the following figures. As displayed in Figure \ref{fig:errvsN}, the mean of the identification errors using dual IOP exhibits asymptotic convergence as data length $N$ increases. More specifically, the box plots of error distributions for dual IOP compared with other benchmark methods are given in Figure \ref{fig:errcomp}. A data length of $N=2^{14}-1$ is selected for illustration. Results indicate that given the same PRBS order and FIR order, dual IOP performs the best as it shows an obvious decrease in the error median compared to others. The medians from dual IOP, dual SLP and dual YP with $\mathbf{G}_a$ are $0.943, 1.01, 1.15$ respectively, meaning dual IOP has a reduction in the median by 7\% and 18\%. Note the performance of dual YP strongly depends on the choice of the initial nominal plant. In some cases, especially when the prior knowledge of the real plant is limited, dual YP may lead to less convincing estimates with large uncertainties (for instance, when $\mathbf{G}_X=\mathbf{G}_c$).

Another important observation is that dual IOP demonstrates an increased precision than dual SLP, as the distribution of $\mathrm{err}(\hat{\mathbf{G}})$ for dual SLP results in a significantly larger variance. This is because dual SLP implicitly increases the dimension of the problem compared to dual IOP. In particular, consider the SISO system in (\ref{eq:emp_sys}), the controller $\mathbf{K}_0$ has $n=2$ internal states. This means the dual system-level parameters $\mathbf{R}_K, \mathbf{N}_K, \mathbf{M}_K$ are MIMO transfer functions, while the dual input-output parameters $\mathbf{W}, \mathbf{X}, \mathbf{Y}, \mathbf{Z}$ are still SISO functions. This means dual SLP always optimizes more variables than dual IOP given the same model order. Therefore, the region between the error quartiles (and also extrema) of dual IOP is narrower than those of dual SLP. 

\begin{figure}
  \centering
  \begin{minipage}[b]{0.5\textwidth}
    \includegraphics[width = \textwidth]{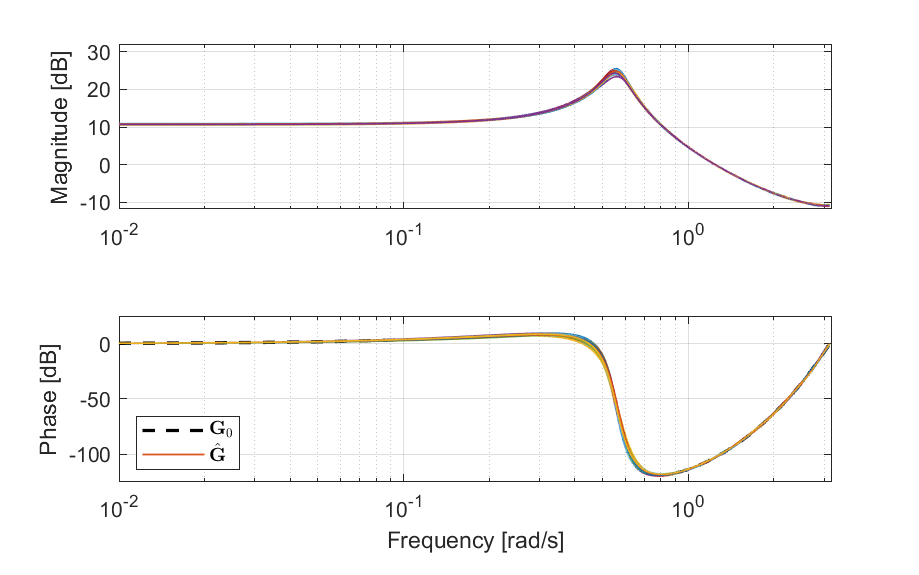}
\caption{Bode plots of 100 independent dual IOP identification trials. All the dual IOP estimates (solid curves) provide very good approximations of the real plant (dashed curve).} 
\label{fig:iop100}
  \end{minipage}
  \hfill
  \begin{minipage}[b]{0.5\textwidth}
    \includegraphics[width = \textwidth]{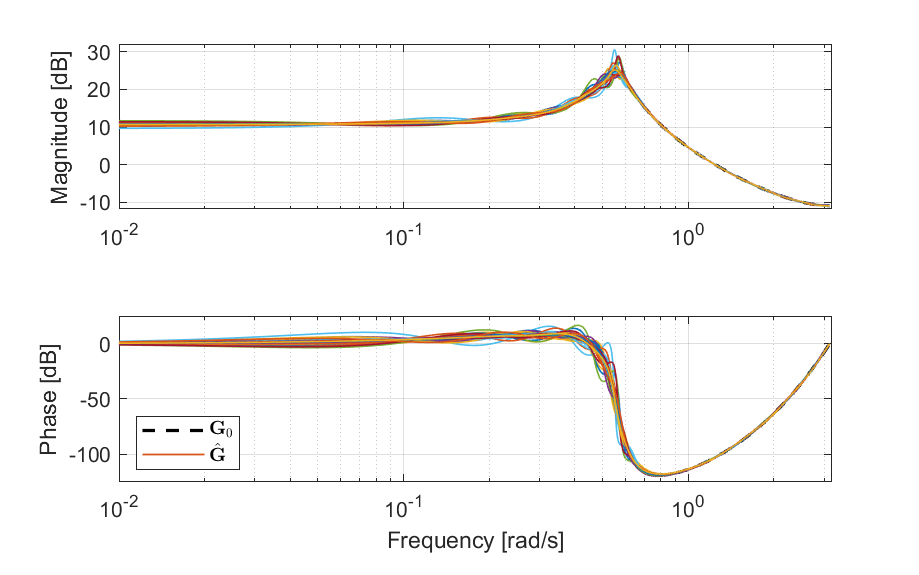}
\caption{Bode plots of 100 independent dual SLP identification trials. The dual SLP estimates (solid curves) can deviate from the true plant for some trials.}
\label{fig:slp100}
  \end{minipage}
\end{figure}

We further compare the plants estimated by dual IOP and dual SLP via Bode plots. The Bode plots of 100 independent experiments using dual IOP and dual SLP are given in Figures \ref{fig:iop100} and \ref{fig:slp100}, respectively. As is evident from the figures, over 100 trials all the estimates by dual IOP are located around the true plant with smaller offsets. On the other hand, dual SLP estimates show larger deviations to the original plant, thus indicating dual SLP is more vulnerable to overfitting and dual IOP offers superior performance compared to the above two benchmark methods. 

\begin{figure}
    \centering
    \includegraphics[width = 0.49\textwidth]{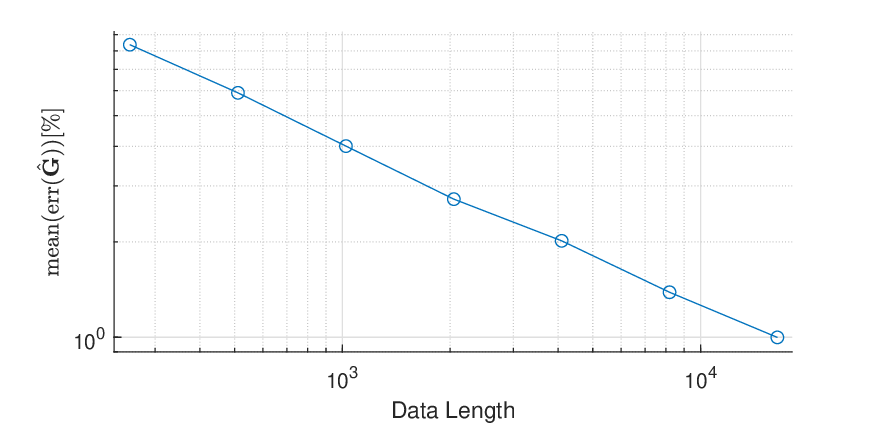}
    \caption{Plots of the mean of $\mathrm{err}(\hat{G})$ versus data length using dual IOP. The error shows asymptotic convergence as the data length increases.}
    \label{fig:errvsN}
\end{figure}

%% file: sections/05_conclu.tex
\section{Conclusion}
\label{sec:conclu}
In this paper, we propose a novel closed-loop identification framework, named dual input-output parameterization (dual IOP), which is the dual problem to the input-output parameterization in control synthesis. Given the knowledge of the controller, the dual IOP identifies the plant by optimizing the closed-loop transfer functions subject to a group of linear equality constraints. The estimated plant is guaranteed to be stabilized by the known controller, and the estimation error shows asymptotic convergence with respect to the input data length. The dual IOP improves upon the benchmark methods since (a) it exhibits decreased empirical mean and variance of the estimation error, and (b) it does not rely on any pre-computations, such as doubly-coprime factorization or state-space realization of the controller. For future work, a possible direction is to investigate the performance of the dual IOP framework in the identification of systems with specific closed-loop behaviors, such as plants that are positively stabilized by the known controller.